\documentclass[10pt,fleqn]{article}
\usepackage{mathptmx}

\input{epsf.tex}
\usepackage{latexsym,amsfonts,amssymb,epsfig,verbatim}
\usepackage{amsmath,amsthm,amssymb,latexsym,graphics,textcomp}
\usepackage{eucal,eufrak}
\usepackage{graphicx}
\usepackage{url}

\theoremstyle{plain}
\newtheorem{theorem}{Theorem}[section]
\newtheorem{lemma}[theorem]{Lemma}
\newtheorem{corollary}[theorem]{Corollary}

\newtheorem{proposition}[theorem]{Proposition}

\newtheorem*{claim}{Claim}

\theoremstyle{definition}
\newtheorem*{remark}{Remark}

\def\1{\mathbf 1}

\def\EL{\mathcal{EL}}
\def\ML{\mathcal{ML}}

\def\PML{\mathbb P \mathcal{ML}}

\def\T{\mathrm{Teich}}
\def\Mod{\mathrm{Mod}}

\def\seg{{\, \mathrm{ seg}}}
\def\ray{{\, \mathrm{ ray}}}

\def\A{\mathcal A}
\def\C{\mathcal C}

\def\T{\mathcal T}

\def\Ca{\mathfrak C}
\def\G{\mathfrak G}
\def\Q{\mathcal Q}

\def\diam{\mathrm{diam}}

\def\d{\mathrm{d}}

\def\Dom{\mathrm{Dom}}

\def\g{\mathcal G}

\begin{document}

\title{\textbf{Uniform convergence in the mapping class group}}
\author{Richard P. Kent IV\thanks{Partially supported by a Donald D. Harrington Dissertation Fellowship and an NSF Postdoctoral Fellowship.} \, and  Christopher J. Leininger\thanks{Partially supported by NSF Grant DMS-0603881}}
\maketitle
{\abstract We characterize convex cocompact subgroups of the
mapping class
group of a surface in terms of uniform convergence actions on the zero
locus of the limit set. We also construct subgroups that act as uniform
convergence groups on their limit sets, but are not convex cocompact. }

\section{Introduction}\label{introsect}

The notion of convex cocompactness for Kleinian groups was extended to subgroups of $\Mod(S)$, the mapping class
group of a closed surface $S$, by Farb and Mosher \cite{FMcc} by way of analogy.  This analogy was strengthened by the
authors in \cite{kentleininger} (see also Hamenst\"adt \cite{hamenstadt}), appealing to the work of McCarthy and
Papadopoulos in \cite{mccarthypapa}, which itself describes an analogy between the dynamical aspect of Kleinian groups
and that of subgroups of $\Mod(S)$.

In conversation at the 2005 Ahlfors--Bers colloquium, Ed Taylor asked us whether there is a formulation of convex
cocompactness for mapping class groups, analogous to the following for Kleinian groups; see e.g. \cite{freden}.

\begin{theorem} \label{kleintrip}
A nonelementary Kleinian group $\Gamma$ is convex cocompact if and only if the action of $\Gamma$ on the limit set
$\Lambda_\Gamma$ is a uniform convergence action.
\end{theorem}

Recall that an action of a group $G$ on a perfect compact metrizable space $X$ is a (discrete) {\em convergence action}
if the diagonal action on the space of distinct triples in $X$ is properly discontinuous, and that it is {\em uniform}
if this associated action is cocompact; see e.g. \cite{GeMar}, \cite{tukia}, and \cite{bowditchcon}. Bowditch has shown
that uniform convergence groups provide an alternative characterization of word-hyperbolicity for a group
\cite{bowditchtopchar}.

Although certain aspects of the theory of Kleinian groups have very nice analogies in the world of mapping class groups, there are limitations to this; see e.g. \cite{masurclass}, \cite{masurwolf}, \cite{brockfarb}, \cite{AAS05}, \cite{BDM05}.
Our first theorem describes another such limitation.

\begin{theorem} \label{toobad}
There exist nonelementary irreducible subgroups $G < \Mod(S)$ which act as uniform convergence groups on their limit
set $\Lambda_G$, but are not convex cocompact. In fact, there exists a rank--two free subgroup $G < \Mod(S)$ which is
not convex cocompact, but which admits a $G$--equivariant parametrization of its limit set
\[ \partial G \to \Lambda_G. \]
\end{theorem}

When presented a property of Kleinian groups involving the limit set, an analogous property of subgroups of $\Mod(S)$
will often involve the zero locus of the limit set $Z \Lambda_G$ in its place. This is of a certain necessity, illustrated by the
theorem above; And a certain appropriateness, illustrated by the theorem below.

\begin{theorem} \label{tripleZ}
Let $G$ be a nonelementary subgroup of $\Mod(S)$.  Then $G$ is convex cocompact if and only if $G$ acts as a uniform
convergence group on $Z \Lambda_G$.
\end{theorem}

Our proof goes through another characterization of convex cocompactness related to the uniform convergence criterion.
This characterization makes use of an alternative space of triples.  To describe it, recall that a pair
$[\lambda_-],[\lambda_+] \in \PML(S)$ is said to \textit{bind} $S$ if $i(\lambda_-,\mu) + i(\lambda_+,\mu) \neq 0$ for
every $\mu \in \ML(S)$.  Now set
\[
\Omega = \Big\{ \big( [\lambda_-],[\lambda_+],[\mu] \big) \in \PML(S)^3 \, \Big| \, [\lambda_-],[\lambda_+] \mbox{ bind
} S \, , i(\lambda_-,\mu) \neq 0 \, \mbox{ and } \, i(\lambda_+,\mu) \neq 0 \Big\}.
\]

For $G < \Mod(S)$, we set $\Omega_G = \Omega \cap \Lambda_G^3$ and prove the following.

\begin{theorem} \label{triplespecial}
Suppose $G < \Mod(S)$ is a nonelementary group. Then $G$ is convex cocompact if and only if $\Omega_G \neq \emptyset$
and $G$ acts cocompactly on it.
\end{theorem}

\subsection{Outline of the paper}

To illustrate the naturality of the space $\Omega_G$, we prove that it is usually nonempty and that $G$ always acts
properly discontinuously on it [Theorem \ref{pd}]. The proof of this latter fact is similar in spirit to the proof of
proper discontinuity for any Kleinian group on the space of distinct triples of points in the limit set. In Section
\ref{tripcharsection} we prove Theorem \ref{triplespecial}. The proof relies on the notion of conical limit point as
defined in \cite{kentleininger} and the corresponding characterization of convex cocompactness. With Theorem
\ref{triplespecial} at our disposal, we prove Theorem \ref{tripleZ} in Section \ref{tripchar2section} appealing to
length/intersection number arguments.

Finally, in Section \ref{countersection}, we construct the counterexamples proving Theorem \ref{toobad}. The examples
given are subgroups generated by a pseudo-Anosov mapping class and a reducible mapping class which is pseudo-Anosov on
the complement of a simple closed curve.

\subsection{Notation and background}

We adhere to the notation, conventions, and definitions of \cite{kentleininger}, with a few minor exceptions that we
spell out here.  For the convenience of the reader, we will recall the notation necessary for the discussion, referring
to \cite{kentleininger} for a more thorough treatment of the background and for the appropriate references.

We write $\T(S)$ for the Teichm\"uller space of $S$, $\ML(S)$ for the space of measured laminations on $S$ and
$\PML(S)$ for its projectivization.  For a subgroup $G < \Mod(S)$, $\Lambda_G$ and $Z \Lambda_G$ denote the limit set
of $G$ and its zero locus, respectively.  A $\pi_1$--injective subsurface $Y \subset S$ is called a domain.  We write
$\xi(Y) = 3g +p$ to denote the complexity.  We will assume throughout that $\xi(S) \geq 5$, since the other cases are
not of interest here.  The complex of curves is denoted $\C(Y)$ and the complex of arcs by $\A(Y)$, with the boundary
at infinity (of both) being $\EL(S)$, the space of ending (filling) laminations.  We will often refer to a vertex (or inappropriately,
a point) of $\C(S)$ as a curve, and vice versa.  We write $\pi_Y$ for the subsurface projection and $\d_Y$ for the
subsurface distance. For a pseudo-Anosov element $f \in \Mod(S)$, we let $[\lambda_+(f)],[\lambda_-(f)]$ denote the
stable and unstable projective measured laminations of $f$, respectively.

In what follows, unlike in \cite{kentleininger}, we do not assume that a uniquely ergodic lamination is filling. Also,
if $v \in \C(S)$ is a curve, then we write $\A(v)$, $\pi_v$, and $\d_v$ in place of $\A(Y)$, $\pi_Y$, and $\d_Y$, were
$Y$ is the annulus with core $v$.\\

\noindent {\bf Acknowledgements.}  The authors would like to thank Ed Taylor for asking us about the relationship with
convergence actions.

\begin{remark}
Fenley and Mosher \cite{fenleymosher} have also studied convex cocompactness in connection with convergence and
uniform convergence actions, but following a different line of questions.  Their work relates properties of actions of $G$
to actions of the associated surface group extension $\Gamma_G$.
\end{remark}

\section{Proper discontinuity on $\Omega$} \label{pdsection}

To motivate the use of $\Omega$ as a replacement for the space of distinct triples, let us prove the following easy fact.

\begin{proposition} \label{nopd}
The action of $\Mod(S)$ on the set of distinct $n$-tuples of points in $\PML(S)$ is not properly discontinuous for any
$n > 0$.
\end{proposition}

\begin{proof}
Let $T$ be a Dehn twist in any simple closed curve $\alpha$ in $S$. There is an uncountable set of points in $\PML(S)$ fixed
pointwise by $T$.  Indeed, there is a positive dimensional subspace of $\PML(S)$ consisting of laminations having zero intersection
number with $\alpha$, and this entire set is fixed by $T$.
Any set of $n$ distinct points in this set determines a point in the space of distinct $n$-tuples fixed by $T$.
This is a compact set fixed by an infinite order element.
\end{proof}

Similar phenomena occur for many subgroups $G < \Mod(S)$ acting on the space of distinct $n$-tuples of points in
$\Lambda_G$.  The spaces $\Omega$ and $\Omega_G$ circumvent this problem.
In contrast to Proposition \ref{nopd}, we have the following.

\begin{theorem} \label{pd}
$\Mod(S)$ acts properly discontinuously on $\Omega$.
\end{theorem}

This immediately implies

\begin{corollary}
$G$ acts properly discontinuously on $\Omega_G$.
\end{corollary}

To prove Theorem \ref{pd}, we wish to construct a $\Mod(S)$--equivariant continuous map $\Pi:\Omega \to \T(S)$.  This
will readily imply that the action of $\Mod(S)$ on $\Omega$ is properly discontinuous, since the action on $\T(S)$ is.

We consider the space of all binding pairs of measured laminations on $S$ with intersection number one
\[
\Q^1(S) = \big\{ (\lambda_-,\lambda_+) \in \ML(S)^2 \, \big| \, \lambda_-,\lambda_+ \mbox{ binds } S \mbox{ and }
i(\lambda_-,\lambda_+) = 1 \big\}.
\]
There is a canonical homeomorphism between this space and the bundle over $\T(S)$ whose fiber at $X$ is the space of
unit norm quadratic differentials, holomorphic with respect to $X$---see \cite{gardinermasur}.  The homeomorphism of
the former with the latter is defined by sending a pair $(\lambda_-,\lambda_+)$ to the unique quadratic differential
having horizontal and vertical foliations naturally associated to the measured laminations $\lambda_-$ and $\lambda_+$,
respectively. We use this homeomorphism to identify these two spaces, no longer distinguishing between the two.

We recall Masur's description of the {\em Teichm\"uller geodesic flow} on $\Q^1(S)$---see \cite{masurtransitivity} and
\cite{masurinterval}. This flow $\varphi_t$ is given by
\[
\varphi_t(\lambda_-,\lambda_+) = (e^{-t} \lambda_-,e^t \lambda_+).
\]
Under the natural projection from $\Q^1(S)$ to $\T(S)$, the flow lines project to geodesics, with $t$ a unit speed
parameter.  Every geodesic arises in this way.

It follows that the space of geodesics in $\T(S)$ (or flow lines on $\Q^1(S)$) is canonically identified with the space
\[
\G = \Big\{ \big( [\lambda_-],[\lambda_+] \big) \in \PML(S)^2 \, \Big| \, [\lambda_-],[\lambda_+] \mbox{ bind } S
\Big\}.
\]

We now describe a map $\widehat\Pi:\Omega \to \Q^1(S)$.  The desired map $\Pi:\Omega \to \T(S)$ is then obtained by
composing with the $\Mod(S)$--equivariant continuous projection $\Q^1(S) \to \T(S)$. For any triple $\big(
[\lambda_-],[\lambda_+],[\mu] \big)$, we consider the flow line $\hat \tau_{[\lambda_-],[\lambda_+]}$ defined by $\big(
[\lambda_-],[\lambda_+] \big)$ and define $\widehat \Pi \big( [\lambda_-],[\lambda_+],[\mu] \big)$ to be the
\emph{balance point} for $[\mu]$ on $\hat \tau_{[\lambda_-],[\lambda_+]}$:  picking representatives
$\lambda_-,\lambda_+ ,\mu$ for the projective classes for which $i(\lambda_-,\lambda_+) = 1$, this is the unique point
$(e^{-t}\lambda_-,e^t\lambda_+) \in \hat \tau_{[\lambda_-],[\lambda_+]}$ for which
\[
i(e^{-t}\lambda_-,\mu) = i(e^t \lambda_+,\mu).
\]
This is independent of choice of representatives.  See \cite{MM1} and \cite{rafi} for more on the notion of balance
point.

The map $\Pi$ is naturally defined, and it is easy to see that it is $\Mod(S)$--equivariant.  Continuity follows
readily from continuity of $i$, but we give the argument for completeness.

\begin{lemma}
$\Pi:\Omega \to \T(S)$ is continuous.
\end{lemma}

\begin{proof}
It suffices to show that $\widehat \Pi$ is continuous. Let $\sigma:\PML(S) \to \ML(S)$ be a continuous section of the
projectivization.  That is, $\sigma$ is continuous and $\sigma [\lambda]$ is a representative of the projective class
of $[\lambda]$. It follows that the map $\hat \sigma: \G \times {\mathbb R} \to \Q^1(S)$ defined by
\[
\hat \sigma \Big( \big( [\lambda_-],[\lambda_+] \big),t \Big) = \left( e^{-t}
\sigma[\lambda_-],\frac{e^t}{i(\sigma[\lambda_-],\sigma[\lambda_+])} \sigma[\lambda_+] \right)
\]
is continuous.

We also consider the map $F: \Omega \times {\mathbb R} \to {\mathbb R}$ defined by
\[
F \Big( \big( [\lambda_-],[\lambda_+],[\mu] \big),t \Big) = \max \left\{ i\left(e^{-t} \sigma[\lambda_-],\sigma[\mu]
\right) \, , \, i \left( \frac{e^t}{i(\sigma[\lambda_-],\sigma[\lambda_+])} \sigma[\lambda_+],\sigma[\mu] \right)
\right\}.
\]
This is continuous, and according to the definition of $\Omega$, is strictly convex as a function of $t$ for every
fixed $\big( [\lambda_-],[\lambda_+],[\mu] \big) \in \Omega$. Therefore, if we set $T \big(
[\lambda_-],[\lambda_+],[\mu] \big)$ to be the unique real number for which $F \Big( \big(
[\lambda_-],[\lambda_+],[\mu] \big),T \big( [\lambda_-],[\lambda_+],[\mu] \big) \Big)$ minimizes the function $F
\Big|_{ \big\{ \big( [\lambda_-],[\lambda_+],[\mu] \big) \big\} \times {\mathbb R}}$, it follows that $T: \Omega \to
{\mathbb R}$ is continuous.

Tracing through the definitions, we see that
\[
\widehat \Pi \big( [\lambda_-],[\lambda_+],[\mu] \big) = \hat \sigma \Big( \big( [\lambda_-],[\lambda_+] \big),T \big(
[\lambda_-],[\lambda_+],[\mu] \big) \Big),
\]
which being composed of continuous functions, is continuous.
\end{proof}

We can now easily prove Theorem \ref{pd}.

\begin{proof}[Proof of Theorem \ref{pd}]
If $K$ is any compact set in $\Omega$, then $\Pi(K)$ is compact. Since the action of $\Mod(S)$ on $\T(S)$ is properly
discontinuous, there are only finitely many elements $g \in \Mod(S)$ for which $g \Pi(K) \cap \Pi(K) \neq \emptyset$.
Since $g\Pi(K) = \Pi(gK)$, and $\Pi(gK) \cap \Pi(K) = \emptyset$ implies $gK \cap K = \emptyset$, it follows that there
are only finitely many $g$ for which $gK \cap K \neq \emptyset$.
\end{proof}

\section{Convex cocompactness I: the action on $\Omega$} \label{tripcharsection}

The goal of this section is to prove Theorem \ref{triplespecial}.  One direction follows from known facts about
hyperbolic groups acting on their boundaries, and the work of Farb and Mosher \cite{FMcc}.  To prove that cocompactness of
the $G$ action on $\Omega_G$ implies convex cocompactness, we will prove that every limit point is conical in the sense
of \cite{kentleininger}. The next lemma is the first ingredient.

\begin{lemma} \label{notemptybp}
The set $\Omega_G \neq \emptyset$ if and only if $G$ is irreducible and nonelementary.  In this case for every $[\lambda] \in
\Lambda_G$ there exists a $[\mu] \in \Lambda_G$, such that $[\lambda],[\mu]$ is a binding pair.
\end{lemma}

\begin{proof}
Suppose $\Omega_G \neq \emptyset$.  Then $G$ is not finite since this implies $\Lambda_G = \emptyset$.  Furthermore,
$G$ cannot be reducible. If it were, then every element of $\Lambda_G$ would have zero intersection number with the
essential reduction system for $G$ (see \cite{mccarthypapa}, Section 7), and hence no pair in $\Lambda_G$ could be
binding.  This is a contradiction.

Conversely, suppose $G$ is irreducible and nonelementary.  By Ivanov's Theorem \cite{ivanov} any irreducible subgroup
contains a pseudo-Anosov element. Because $G$ is nonelementary, there are two pseudo-Anosov elements $g,h \in G$ with
no fixed points in common and so $\big( [\lambda_-(g)],[\lambda_+(g)],[\lambda_-(h)] \big) \in \Omega_G$, proving that
it is nonempty. Moreover, any $[\lambda] \in \Lambda_G$ different than $[\lambda_-(g)]$ binds with $[\lambda_-(g)]$,
and if $[\lambda] = [\lambda_-(g)]$, then $[\lambda],[\lambda_+(g)]$ is a binding pair.
\end{proof}

\begin{proof}[Proof of Theorem \ref{triplespecial}]
If $G$ is convex cocompact, then by Theorem 1.1 of \cite{FMcc}, every lamination in $\Lambda_G$ is filling and uniquely
ergodic and there is a $G$--equivariant homeomorphism $\partial G \to \Lambda_G$.  Therefore $\Omega_G$ is
$G$--equivariantly homeomorphic to the space of distinct triples in $\Lambda_G$.  The action of $G$ on the latter space
is known to be properly discontinuous and cocompact, see \cite{freden} and \cite{bowditchtopchar}, and so the action
on the former is also.

Now suppose $\Omega_G \neq \emptyset$ and the action is cocompact.  Fix $[\lambda_+] \in \Lambda_G$. By Lemma
\ref{notemptybp}, there exists $[\lambda_-] \in \Lambda_G$ so that $[\lambda_-],[\lambda_+]$ is a binding pair. We
choose representatives so that $i(\lambda_-,\lambda_+) = 1$; note that this specifies a parametrization of the geodesic
$\tau_{\lambda_-,\lambda_+}$. Furthermore, since $G$ is irreducible, the set of stable laminations of pseudo-Anosov
elements in $G$ is dense in $\Lambda_G$. Let $\big\{ [\mu(n)] \big\}_{n=1}^{\infty}$ be any sequence of stable
laminations in $\Lambda_G$ converging to $[\lambda_+]$. We choose representatives $\mu(n)$ so that $\mu(n) \to
\lambda_+$ as $n \to \infty$. We may assume that $[\lambda_-] \neq [\mu(n)]$ and $[\lambda_+] \neq [\mu(n)]$ for all
$n$.

It follows that $\Big\{ \big( [\lambda_-],[\lambda_+],[\mu(n)] \big) \Big\}_{n=1}^{\infty} \subset \Omega_G$. Moreover,
this sequence must be diverging in $\Omega_G$ since its limit in $\Lambda_G^3$ is the point $\big(
[\lambda_-],[\lambda_+],[\lambda_+] \big) \not \in \Omega_G$. Therefore, by cocompactness of the $G$--action on
$\Omega_G$, there exists a sequence $g_n \in G$ so that
\[
g_n \big( [\lambda_-],[\lambda_+],[\mu(n)] \big) = \big( g_n[\lambda_-],g_n[\lambda_+],g_n[\mu(n)] \big) \to \big(
[\lambda_-(\infty)],[\lambda_+(\infty)],[\mu(\infty)] \big)
\]
as $n$ tends to infinity. Since $\Big\{ \big( [\lambda_-],[\lambda_+],[\mu(n)] \big) \Big\}$ diverges, we may assume,
by passing to a subsequence if necessary, that the $g_n$ are all distinct.

Since $\Pi$ is continuous, we see that as $n \to \infty$
\[
\Pi \Big( g_n \big( [\lambda_-],[\lambda_+],[\mu(n)] \big) \Big) \to \Pi \big(
[\lambda_-(\infty)],[\lambda_+(\infty)],[\mu(\infty)] \big).
\]
By passing to a further subsequence, we may assume that
\[
\d_\T \Big( \Pi \Big( [\lambda_-(\infty)],[\lambda_+(\infty)],[\mu(\infty)] \Big) ,\Pi \Big( g_n \big(
[\lambda_-],[\lambda_+],[\mu(n)] \big) \Big) \Big) \leq 1
\]
for every $n$.

Since $G$ acts by isometries on $\T(S)$, and since $\Pi$ is $G$--equivariant this says that
\begin{equation} \label{closetoray}
\d_\T \Big( g_n^{-1} \Big( \Pi \big( [\lambda_-(\infty)],[\lambda_+(\infty)],[\mu(\infty)] \big) \Big),\Pi \big(
[\lambda_-],[\lambda_+],[\mu(n)] \big) \Big) \leq 1.
\end{equation}

Now consider the ray
\[
\vec{\tau}_{\lambda_-,\lambda_+} = \big\{ \tau_{\lambda_-,\lambda_+}(t) \, | \, t \geq 0  \big\}.
\]
If we can show that $\Pi \big( [\lambda_-],[\lambda_+],[\mu(n)] \big)$ is contained in this ray for sufficiently large
$n$, then (\ref{closetoray}) implies that the tail of the sequence
\[
\left\{ g_n^{-1} \Big( \Pi \big( [\lambda_-{\infty}],[\lambda_+(\infty)],[\mu(\infty)] \big) \Big) \right\}
\]
provides infinitely many points of the $G$--orbit of $\Pi \big( [\lambda_-(\infty)],[\lambda_+(\infty)],[\mu(\infty)]
\big)$ within a distance $1$ of this ray. Since the direction of $\vec{\tau}_{\lambda_-,\lambda_+}$ is $[\lambda_+]$,
this will show that $[\lambda_+]$ is a conical limit point, and the proof will be complete.

By definition, $\Pi \big( [\lambda_-],[\lambda_+],[\mu(n)] \big)$ are all points on $\tau_{\lambda_-,\lambda_+}$. As $n
\to \infty$, we have $i(\lambda_+,\mu(n)) \to 0$ and eventually $i(\lambda_-,\mu(n)) \geq \frac{1}{2}$. This follows
from continuity of $i$, the fact that $\mu(n) \to \lambda_+$, and $i(\lambda_-, \lambda_+) = 1$. In particular, this
implies that for sufficiently large $n$, we have
\[
i(\lambda_+,\mu(n)) < i(\lambda_-,\mu(n)).
\]
Since $\Pi \big( [\lambda_-],[\lambda_+],[\mu(n)] \big)$ is the point $\tau_{\lambda_-,\lambda_+}(t)$ for which
\[
e^t i(\lambda_+,\mu(n)) = e^{-t} i(\lambda_-,\mu(n))
\]
we see that for all sufficiently large $n$, we must have $t > 0$ and hence $\Pi \big( [\lambda_-],[\lambda_+],[\mu(n)]
\big)$ lies on $\vec{\tau}_{\lambda_-,\lambda_+}$ as required.
\end{proof}

\section{Convex cocompactness II: the zero locus} \label{tripchar2section}

We will need an alternate description of a uniform convergence action---see \cite{bowditchtopchar}.

\begin{theorem}[Bowditch,Tukkia] \label{alternatedescription} The action of a group $G$ on $X$ is a convergence
action if and only if for every sequence $\{g_n\}$ of distinct elements of $G$ there is a subsequence $\{g_{n_k}\}$ and
a point $x \in X$ so that the restriction of $g_{n_k}$ to $X - \{x\}$ converges uniformly on compact sets to a constant
function.

A convergence action of $G$ on $X$ is uniform if and only if for every $x \in X$, there exists a pair of distinct
points $a,b \in X$ and a sequence $\{g_n\}$ so that $\lim g_n(x) = a$ and the restriction of $g_n$ to $X - \{x\}$
converges uniformly on compact sets to the constant function with value $b$.
\end{theorem}

Using this, we now prove
\begin{lemma} \label{edgeslemma}
If $G$ acts as a uniform convergence group on $Z \Lambda_G$, then every lamination in $\Lambda_G$ is filling and uniquely ergodic.
In particular, $Z \Lambda_G = \Lambda_G \neq \PML(S)$.
\end{lemma}

\begin{proof}  If the conclusion of the lemma were false, then $Z \Lambda_G$ would contain a positive dimensional projective
simplex of measures. Consider an edge of this simplex given by $\big\{ [t\mu + (1-t)\lambda] \, \big| \, t \in [0,1]
\big\}$ for some $[\lambda],[\mu] \in Z \Lambda_G$ with $i(\lambda,\mu) = 0$.

Proposition \ref{alternatedescription} implies that there is a sequence $\{g_n \} \subset G$ so that as $n \to \infty$
\begin{equation} \label{conveqn1}
g_n \left[ \frac{\mu}{2} + \frac{\lambda}{2} \right] \to [\eta_1] \in Z \Lambda_G
\end{equation}
and $g_n$ converges uniformly on compact sets to the constant map with value \linebreak $[\eta_2] \in Z \Lambda_G -
\big\{ [\eta_1] \big\}$:
\begin{equation} \label{conveqn2}
g_n \Big|_{Z \Lambda_G - \left\{ \left[ \frac{\mu}{2} + \frac{\lambda}{2} \right] \right\} } \to [\eta_2] \neq
[\eta_1].
\end{equation}

Fix a hyperbolic metric $X$ on $S$ and let $\{t_n\}$ be positive numbers so that
\[1 = t_n \ell_X\left( g_n \left( \frac{\mu}{2}  + \frac{\lambda}{2} \right) \right) = \ell_X \left( \frac{t_n g_n \mu}{2} \right) + \ell_X \left( \frac{t_n g_n \lambda}{2} \right) \]
It follows that both of the lengths $\ell_X(t_n g_n \mu)$ and $\ell_X(t_n g_n \lambda)$ must be bounded above, and
after passing to a subsequence, at least one of them is bounded below by a positive number.  We can therefore pass to a
subsequence so that at least one of $\{t_n g_n \mu\}$ and $\{t_n g_n \lambda\}$ converges,
and that if only one of these sequence converges, then the length in $X$ of the other must tend to zero.\\

\noindent \textbf{Case 1.}  After passing to a subsequence both $\{t_n g_n \mu\}$ and $\{t_n g_n \lambda\}$
converge to laminations $\eta_3$ and $\eta_4$ in $\ML(S)$.\\

According to (\ref{conveqn2}), $[\eta_3] = [\eta_2] = [\eta_4]$.  But then $[\eta_3 + \eta_4] = [\eta_2]$ and combining
this with (\ref{conveqn1}) we have
\begin{eqnarray*}
[\eta_1] & = & \lim_{n \to \infty} g_n \left[ \frac{\mu}{2} + \frac{\lambda}{2} \right]\\
 & = & \lim_{n \to \infty} g_n
\left[ t_n \left( \frac{\mu}{2} + \frac{\lambda}{2} \right) \right]\\
 & = & \lim_{n \to \infty} \left[\frac{t_n g_n \mu}{2} + \frac{t_n g_n \lambda}{2} \right]\\
  & = & [\eta_3 + \eta_4]\\
  & = & [\eta_2].
\end{eqnarray*}
This is a contradiction since $[\eta_1] \neq [\eta_2]$.\\

\noindent \textbf{Case 2.}  After passing to a subsequence, only one of the sequences, $\{t_n g_n \mu\}$, say,
converges to a lamination $\eta_3$ in
$\ML(S)$, and $\ell_X(t_n g_n \lambda) \to 0$ as $n \to \infty$.\\

According to (\ref{conveqn2}) we must have $[\eta_3] = [\eta_2]$.  Since $\ell_X(t_n g_n \lambda) \to 0$, we see that
\[
[\eta_1] = \lim_{n \to \infty} \left[ \frac{t_n g_n \mu}{2} + \frac{t_n g_n \lambda}{2} \right] = \lim_{n \to \infty}
\left[\frac{t_n g_n \mu}{2} \right] = \left[ \frac{\eta_3}{2} \right] = [\eta_2].
\]
This is also contradicts the fact that $[\eta_1] \neq [\eta_2]$.\\

These two cases exhaust the possibilities, so $Z \Lambda_G$ can contain no positive dimensional simplices and hence all
laminations in $\Lambda_G$ are filling and uniquely ergodic.
\end{proof}

Theorem \ref{tripleZ} now follows from this lemma and Theorem \ref{triplespecial}.

\begin{proof}[Proof of Theorem \ref{tripleZ}.]  Suppose first that $G$ is convex cocompact.  Since $Z \Lambda_G = \Lambda_G$, as in the proof of Theorem \ref{triplespecial}, the theorem follows from \cite{FMcc} and the fact that the action of a hyperbolic group on its Gromov boundary is a uniform convergence action.

To prove the other direction, suppose $G$ acts as a uniform convergence group on $Z \Lambda_G$.
According to Lemma \ref{edgeslemma} every limit point is filling and uniquely ergodic, and so $Z \Lambda_G = \Lambda_G$ and the space of distinct triples is equal to $\Omega_G$.
Moreover, because $G$ is nonelementary, $\Omega_G$ is nonempty.
Therefore, $G$ acts cocompactly on $\Omega_G \neq \emptyset$, and so by Theorem \ref{tripleZ}, $G$ is convex cocompact.
\end{proof}

\section{Examples} \label{countersection}

In this section we describe the construction of the examples proving Theorem \ref{toobad}.
We begin with a few technical results on geodesics in the curve complex and their ending laminations.

For consistency, we follow the convention that the distance between two sets in the curve complex, like the distance
between the subsurface projections, is the diameter of their union.  We write $Y$ and $Z$ for domains in $S$.


\subsection{Geodesics in the curve complex} \label{factsection}

The following theorem of \cite{MM2} plays a central role in our construction.

\begin{theorem} [{\bf Masur-Minsky}] \label{bgithrm}
Let $Z \subset S$ be a domain.  There exists a constant $M = M(\xi(Z))$ with the following property. Let $Y$ be a
proper connected subdomain of $Z$ with $\xi(Y) \neq 3$ and let $\g$ be a geodesic segment, ray, or biinfinite line in
$\C(Z)$, such that $\pi_{Y}(v) \neq \emptyset$ for every vertex $v$ of $\g$. Then
\[ \diam_{Y}(\g) \leq M. \]
\end{theorem}

One consequence of this theorem that we will require is the following.  It turns local information about a sequence of
vertices in $\C(S)$ into the global statement that the vertices are distinct vertices of a geodesic.

\begin{proposition} \label{geoconstruct}
Suppose $\{v_i\}$ is sequence of vertices of $\C(S)$ (finite, infinite, or biinfinite) such that each $v_i$ is
nonseparating with $Y_i = S \setminus v_i$, $i(v_i,v_{i+1}) \neq \emptyset$ for all $i$, and
\[\d_{Y_i}(v_{i-1},v_{i+1}) > 3M \]
for all $i$.
Then the path in $\C(S)$ obtained by concatenating geodesic segments $[v_i,v_{i+1}]$ is a geodesic.
\end{proposition}
\begin{proof}  The proposition is easily implied by the following stronger statement.
\begin{claim} For any finite number of consecutive vertices $\{v_i\}_{i=j}^k$, any geodesic from $v_j$ to $v_k$ is a concatenation of geodesic segments $[v_i,v_{i+1}]$ for $i = j,...,k-1$.\end{claim}
\begin{proof}
The proof is by induction on $k-j$, with base case $k = j+2$.
By assumption $\d_{Y_{j+1}}(v_j,v_{j+2}) > M$, and so Theorem \ref{bgithrm} implies that any geodesic $[v_j,v_{j+2}]$ must have some vertex with an empty projection to $Y_{j+1}$.
Since $Y_{j+1}$ is the complement $S \setminus v_{j+1}$ and is connected, this is only possible if $v_{j+1}$ is a vertex of $[v_j,v_{j+2}]$.  That is, the geodesic from $v_j$ to $v_{j+2}$ is the concatenation of geodesic segments $[v_j,v_{j+1}]$ and $[v_{j+1},v_{j+2}]$, as required.

Now suppose the claim holds for $k - j \leq n$ and we prove it for $k - j = n+1$.

Fix any $i$ with $j < i < k$ where $k-j = n+1$.  Let $[v_j,v_i]$ and $[v_i,v_k]$ be any geodesic segments.  It follows from the inductive hypothesis that these can be expressed as concatenations of (some possibly degenerate) geodesic segments
\[
[v_j,v_i] = [v_j,v_{i-1}] \cup [v_{i-1},v_i] \quad \mbox{ and } \quad [v_i,v_k] = [v_i,v_{i+1}] \cup [v_{i+1},v_k].
\]
It follows from Theorem \ref{bgithrm} that
\[
\diam_{Y_i}([v_j,v_{i-1}]) \leq M \quad \mbox{ and } \quad \diam_{Y_i}([v_{i+1},v_k]) \leq M.
\]
From this we see that
\begin{eqnarray*}
\d_{Y_i}(v_j,v_k) & \geq & \d_{Y_i}(v_{i-1},v_{i+1}) - \diam_{Y_i}([v_j,v_{i-1}]) - \diam_{Y_i}([v_{i+1},v_k])\\
 & > & 3M - 2M\\
 & = & M.
\end{eqnarray*}

So by Theorem \ref{bgithrm}, any geodesic from $v_j$ to $v_k$ must contain $v_i$ and is therefore a concatenation of geodesic segments $[v_j,v_i]$ and $[v_i,v_k]$.  By induction each of $[v_j,v_i]$ and $[v_i,v_k]$ are concatenations of the required form, and this proves the claim.
\end{proof}
This completes the proof of the proposition.
\end{proof}

We will also need a means of deciding when a filling lamination is uniquely ergodic. We combine Masur's condition for
unique ergodicity proved in \cite{masurhaus} with work of Rafi \cite{rafi} and Masur-Minsky \cite{MM2} to obtain the
necessary criterion.

Given $\mu,\lambda \in \ML(S)$ and $D > 0$, define a set of proper subdomains of $S$ by
\[
\Dom(\mu,\lambda,D) = \{ Z \subset S \, \, | \, \, \pi_Z(\mu) \neq \emptyset \neq \pi_Z(\lambda) \mbox{ and }
\d_Z(\mu,\lambda) > D \}.
\]
If $v \in \C(S)$, we will also use $v$ to denote the lamination supported on the curve $v$ equipped with the transverse counting measure.

\begin{theorem} \label{uergodic}
Suppose $\mu$ is a filling lamination and $v \in \C(S)$ is such that there exists $D > 0$ so that $\Dom(\mu,v,D)$ can
be partitioned into finite families
\[
\Dom(\mu,v,D) = \bigcup_{i \in {\mathbb Z}} \{Z_\alpha \}_{\alpha \in J_i} \quad \mbox{ with } \quad |J_i| < \infty
\]
with the property that for all $i \neq j$,  all $\alpha \in J_i$ and all $\beta \in J_j$ we have
\[
\d(\partial Z_{\alpha}, \partial Z_{\beta}) \geq 4.
\]
Then $\mu$ is uniquely ergodic.

In fact, any Teichm\"uller geodesic ray defined by a quadratic differential with vertical foliation $\mu$ returns to
the thick part of Teichm\"uller space infinitely often.
\end{theorem}

\noindent
{\em Proof.}  Let $\tau = \tau_{v,\mu}$ be a Teichm\"uller geodesic with horizontal and vertical foliations naturally associated to $v$ and $\mu$,
respectively.  Fixing a point on $\tau$, we obtain two geodesic rays, and we denote the one in the positive direction by $\vec{\tau}$.
In \cite{masurhaus}, Masur proves that if $\vec{\tau}$ returns to the thick part of $\T(S)$ infinitely often, then
$\mu$ is uniquely ergodic.\\

\noindent
{\bf Claim.} There exists $C > 0$ so that if a curve $u$ has length less than $C$ along $\tau$, then
\begin{equation} \label{clustershort}
\d(u,\partial Z_\alpha) \leq 1
\end{equation}
for some $i$ and some $\alpha \in J_i$.\\

Assuming the claim, we prove the theorem.
Thus suppose $\vec{\tau}$ exits every thick part of $\T(S)$.
It follows that there exists a sequence of curves $\{u_n\}$ and a discrete set of points $\vec{\tau}(t_n)$ along $\vec{\tau}$ such that
\begin{itemize}
\item the length of $u_n$ at $\vec{\tau}(t_n)$ is less than $C$,
\item $\d(u_n,u_{n+1}) = 1$,
\item $\vec{\tau}(t_n) \to \infty$ as $n \to \infty$.
\end{itemize}
According to \cite{MM1} the sequence $\{ u_n \}$ lies on a quasi-geodesic ray in $\C(S)$.  Moreover, in
\cite{klarreich}, Klarreich shows that $u_n \to |\mu|$ in $\C(S) \cup \EL(S)$ as $n \to \infty$.  Here $\EL(S)$ is the
space of ending laminations (unmeasured filling laminations with the quotient topology from $\ML(S)$) and Klarreich
shows that it is homeomorphic to the Gromov boundary of $\C(S)$. For every $n$ the claim states that there exists
$i(n)$ and $\alpha(n) \in J_{i(n)}$ so that
\[
d(u_n,\partial Z_{\alpha(n)}) \leq 1.
\]
Therefore
\[
\d(\partial Z_{\alpha(n)},\partial Z_{\alpha(n+1)}) \leq d(u_n,u_{n+1}) + 2 \leq 3
\]
and so by induction on $n$ and the hypothesis of the theorem, for all $n \geq 1$ we have $i(n) = i(1)$.  Finiteness of
$J_{i(1)}$ implies $\{u_n\}$ is a bounded sequence in $\C(S)$, contradicting the fact that it converges to $|\mu|$.\\

\noindent {\em Proof of claim.} The proof is similar to Rafi's proof of Theorem 1.5 in \cite{rafi}. The work of Masur
and Minsky (see \cite{rafi}, Theorem 7.3) implies there exists a $K > 0$ so that if $i_Y(v,\mu) > K$, then
\[ \d_Z(v,\mu) > D\]
for some subdomain $Z \subset Y$.

Rafi's Theorem characterizing short geodesics on $\tau$ \cite{rafi} (along with his footnote on page 198) implies that
there exists a constant $C> 0$ so that if $u$ is a curve with length less than $C$ at some point along $\vec \tau$,
then there exists a component $Y$ of $S \setminus u$ (possibly the annulus with core $u$) so that $i_Y(v,\mu) > K$. It
follows that there is a subdomain $Z \subset Y$ so that
\[ \d_Z(v,\mu) > D. \]
Since $u$ is disjoint from $Z$, it follows that $\d(u,\partial Z) \leq 1$.
By hypothesis, $Z = Z_\alpha$ for some $\alpha \in J_i$ and some $i$, proving the claim. \hfill $\Box$\\

\begin{corollary} \label{itworks}
Suppose $\g$ is a geodesic ray in $\C(S)$ with initial vertex $v_0$ containing a sequence of vertices $v_0 < v_1 < v_2
< v_3 < v_4 < \cdots$ so that
\begin{itemize}
\item $\d(v_{2i-1},v_{2i}) \geq 6$ for all $i \geq 1$, and
\item for some $R > 0 $ and all $i \geq 0$ we have $\d_Z(v_{2i-1},v_{2i}) \leq R$ for every domain $Z$ with $\pi_Z(v_j) \neq \emptyset$ for all $j \geq 0$.
\end{itemize}
Then the limiting lamination $|\mu| \in \EL(S)$ of $\g$ is uniquely ergodic.
\end{corollary}

The second hypothesis says that there are no large projection coefficients $\d_Z(v_{2i-1},v_{2i})$.

\begin{proof}
Set $D = 2M + R$ and suppose that $Z$ is some domain for which
\[ \d_Z(v_0,\mu) > D. \]
By Theorem \ref{bgithrm}, $\g$ must pass through the $1$--neighborhood of $\partial Z$.

Suppose that $\pi_Z(v_j) \neq \emptyset$ for every $j \geq 0$.  Let $w \in [v_{2i-1},v_{2i}]$ for some $i \geq 1$ be any vertex. Then
by the triangle inequality and Theorem \ref{bgithrm}
\[
\d_Z(v_{2i-1},v_{2i}) \geq \d_Z(v_0,\mu) - \big( \diam_Z([v_0,v_{2i-1}]) + \diam_Z([v_{2i},\mu]) \big) > D - 2M = R.
\]
This contradicts the hypothesis of the corollary, so it must be that either $\pi_Z(v_j) = \emptyset$ for some $j \geq
0$ or else $w \not \in [v_{2i-1},v_{2i}]$, for every $i \geq 1$.

It follows that for any $Z$ with $\d_Z(v_0,\mu) > D$, we have $\d(\partial Z,w) \leq 1$ for some $w \in
[v_{2i},v_{2i+1}]$ and some $i \geq 0$.  We can therefore partition the domains $Z$ with $\d_Z(v_0,\mu) > D$ into a countable collection of sets
$\{ P_i \}_{i \geq 0}$, so that $Z \in P_i$ if $\partial Z$ lies in the $1$--neighborhood of $[v_{2i},v_{2i+1}]$.  It follows that if $i < j$, $Z \in P_i$, $Z' \in P_j$ then an application of the triangle inequality implies
\[
\d(\partial Z,\partial Z') \geq \d(v_{2i+1},v_{2j}) - 2 \geq 6 - 2 = 4.
\]

This partition satisfies the hypothesis of Theorem \ref{uergodic}, and so $|\mu|$ is uniquely ergodic.
\end{proof}

\subsection{The groups and their properties} \label{examplesection}

We are now ready to describe the construction.

Fix a nonseparating curve $w \in \C(S)$ and suppose $f,h \in \Mod(S)$ satisfy the following conditions:
\begin{enumerate}
\item \label{hreduce} $h$ is reducible, leaves $w$ invariant, and is pseudo-Anosov on $Y = S \setminus w$.
\item \label{htranslation} The translation distance of $h$ on $\A(Y)$ is greater than $3M$ and $h$ fixes a point of $\A(w)$.
\item \label{fpseudo} $f$ is pseudo-Anosov and leaves a geodesic $\g$ in $\C(S)$ invariant.
\item \label{ftranslation} The translation distance of $f$ on $\C(S)$ is at least $6$.
\item \label{vwdistance} There exists a nonseparating curve $v \in \g$ with $\d(w,v) \geq 2$ and $\d(w,v') \geq \d(w,v)$ for all $v' \in \g$.
\item \label{vwprojections} Setting $Z = S \setminus v$, $\d_Z(w,f^k(v)) > 3M$ for all $k \in {\mathbb Z}$, $k \neq 0$.
\end{enumerate}
Here $M$ denotes the constant from Theorem \ref{bgithrm}.
We defer the proof of the existence of such a pair $f$ and $h$ to Section \ref{theyexist}.\\

Let $\Ca$ denote the Cayley graph of the rank-2 free group $\langle f, h \rangle$ {\em abstractly} generated by $f$ and
$h$. There is a canonical homomorphism from $\langle f,h \rangle$ to $\Mod(S)$, and we refer to the image as $G <
\Mod(S)$. We will denote vertices of $\Ca$ by the elements of $\langle f,h \rangle$ which label them.

Theorem \ref{toobad} will follow from

\begin{theorem} \label{examplethrm}
The canonical homomorphism $\langle f,h \rangle \to \Mod(S)$ is injective and there is an $\langle f,h
\rangle$--equivariant homeomorphism
\[ \partial \Ca \to \Lambda_G \subset \PML. \]
Moreover, every element not conjugate to a power of $h$ is pseudo-Anosov.
\end{theorem}

This clearly implies the second part of Theorem \ref{toobad}. The first part follows from the second since a hyperbolic
group acts as a uniform convergence group on its Gromov boundary---see \cite{freden}, \cite{bowditchtopchar}.

\begin{remark}  It is possible to prove Theorem \ref{examplethrm} with fewer conditions imposed on $f$ and $h$ than we
have listed above.   However, these conditions help to simplify the proof.   It is likely true that given any
pseudo-Anosov $f$ and reducible $h$ which is pseudo-Anosov on a subsurface,
sufficiently large powers $f^n$ and $h^n$ will generate a group satisfying Theorem \ref{examplethrm}, but we do not know how to prove this.\\
\end{remark}

Define a map
\[ \Phi: \Ca \to \C(S) \]
by first defining it on vertices as the orbit map so that $\Phi(\1) = v$.
To define $\Phi$ on edges, note first that the segment of $\g$ connecting $v$ to $f(v)$ is a geodesic we denote $[v,f(v)]$.
Define $\Phi([\1,f])$ to be this segment, and extend the map $\langle f,h \rangle$--equivariantly to the rest of the $f$--edges of $\Ca$.
For the $h$--edges, note that any geodesic from $v$ to $h(v)$ must pass through $w$ by hypothesis (\ref{htranslation}) and Theorem \ref{bgithrm}.
We pick any such geodesic and denote it $[v,w][w,h(v)]$ to emphasize that it is a concatenation of the two geodesic segments $[v,w]$ and $[w,h(v)]$.
We then define $\Phi([\1,h])$ to be this geodesic, and extend to the rest of the $h$--edges of $\Ca$ $\langle f,h \rangle$--equivariantly.

A geodesic in $\Ca$ all of whose edges are $f$--edges will be called an {\em $f$--geodesic}.
Likewise, any geodesic with all $h$--edges will be called an {\em $h$--geodesic}.

We first observe that the $\Phi$--image of any $f$--geodesic is a geodesic as it is simply an $\langle f,h \rangle$--translate of the segment of $\g$
from $v$ to some $f^k(v)$.

On the other hand, the $\Phi$--image of an $h$--geodesic is only a geodesic in the simplest case: when the
$h$--geodesic is an $h$-edge. To see this, note that the geodesic is an $\langle f,h \rangle$--translate of the path
\[
[v,w][w,h(v)][h(v),h(w)][h(w),h^2(v)]\cdots[h^k(v),h^k(w)][h^k(w),h^{k+1}(v)]\]
\[
 = [v,w][w,h(v)][h(v),w][w,h^2(v)]\cdots[h^k(v),w][w,h^{k+1}(v)]
\]
where the equality comes from hypothesis (\ref{hreduce}) that $h(w) = w$.
We can {\em straighten} this to a geodesic segment by simply deleting the middle portion
\[
[w,h(v)][h(v),w][w,h^2(v)]\cdots[h^k(v),w]
\]
from the path.
Note that the result $[v,w][w,h^k(v)]$ is indeed a geodesic, again by hypothesis (\ref{htranslation}) and Theorem \ref{bgithrm}.

We call $v$, $w$, and $h^k(v)$ the {\em special vertices} of $[v,w][w,h^k(v)]$.
The straightening of the $\Phi$--image of any $h$--geodesic has the form $\varphi([v,w][w,h^k(v)])$ for some $\varphi \in \langle f,h \rangle$
and we call the vertices $\varphi(v)$, $\varphi(w)$, and $\varphi h^k(v)$ the {\em special vertices} of this straightening.
We also refer to the endpoints of the $\Phi$--image of any $f$--geodesic as its {\em special vertices}.

Given a geodesic segment $\gamma$ of $\Ca$, we define the {\em straightening} of $\Phi(\gamma)$, denoted $\Phi_*(\gamma)$,  by first writing
it as an alternating concatenation of $\Phi$--images of $f$--geodesics and $h$--geodesics, then straightening each of the $\Phi$--images
of the $h$--geodesics.
Assuming that $\gamma$ starts with an $f$--geodesic, we denote the set of special vertices of
$\Phi_*(\gamma)$ by $\{v_1,v_2,w_3,v_4,v_5,w_6,... \}$.
If $\gamma$ starts with an $h$--geodesic, then we denote the set of special vertices of $\Phi_*(\gamma)$ by $\{v_1,w_2,v_3,v_4,w_5,v_6,... \}$.
Here consecutive vertices $v_i,v_{i+1}$ are the special vertices of the $\Phi$--image of an $f$--geodesic, while consecutive triples
$v_{i-1},w_i,v_{i+1}$ are the special vertices of the straightening of the $\Phi$--image of an $h$--geodesic.

\begin{lemma} \label{straightgeos}
For every geodesic segment $\gamma$ in $\Ca$, $\Phi_*(\gamma)$ is a geodesic in $\C(S)$. If $\gamma$ has positive
length, then so does $\Phi_*(\gamma)$.
\end{lemma}

\begin{proof}
We wish to apply Proposition \ref{geoconstruct} to the sequence of special vertices, so we need to understand the local picture of three consecutive
vertices.
There are three situations to analyze:
\[ \{v_{i-1},v_i,w_{i+1}\} \, , \quad \{v_{i-1},w_i,v_{i+1} \} \, , \quad \{w_{i-1},v_i,v_{i+1} \}. \]
We write the complementary domains for the special vertices as
\[ Y_i = S \setminus w_i \quad \mbox{ and } \quad Z_i = S \setminus v_i\]
which are just $\langle f,h \rangle$--translates of $Y$ and $Z$, respectively.\\

\noindent
{\bf Case 1.} $\{v_{i-1},v_i,w_{i+1} \}$.\\
There is an element $\varphi \in \langle f,h \rangle$ and $k \neq 0$ so that
\[ v_{i-1} = \varphi(v) \, , \quad v_i = \varphi f^k(v) \, , \quad w_{i+1} = \varphi f^k (w).\]
It follows from hypothesis (\ref{ftranslation}) that $\d(v_{i-1},v_i) = \d(v,f^k(v)) \geq 6$, and in particular
\linebreak $i(v_{i-1},v_i) \neq 0$.

Also, since $Z_i = \varphi f^k(Z)$, hypothesis (\ref{vwprojections}) implies
\begin{eqnarray*}
\d_{Z_i}(v_{i-1},w_{i+1}) & = &\d_{\varphi f^k(Z)}(\varphi(v),\varphi f^k(w))\\
 & = & \d_Z(f^{-k}(v),w)\\
 & >& 3M.\\
\end{eqnarray*}

\noindent
{\bf Case 2.} $\{ v_{i-1},w_i,v_{i+1} \}$.\\
There is an element $\varphi \in \langle f,h \rangle$ and $k \neq 0$ so that
\[ v_{i-1} = \varphi(v) \, , \quad w_i = \varphi(w) \, , \quad v_{i+1} = \varphi h^k(v).\]
It follows from hypothesis (\ref{vwdistance}) that $\d(v_{i-1},w_i) = \d(v,w) \geq 2$, and so again $i(v_{i-1},w_i) \neq 0$.

Since $Y_i = \varphi(Y)$, hypothesis (\ref{htranslation}) implies
\begin{eqnarray*}
\d_{Y_i}(v_{i-1},v_{i+1}) & = & \d_{\varphi(Y)}(\varphi(v),\varphi h^k(v))\\
 & = & \d_Y(v,h^k(v))\\
 &> & 3M.\\
\end{eqnarray*}

\noindent
{\bf Case 3.} $\{ w_{i-1},v_i,v_{i+1} \}$.\\
There is an element $\varphi \in \langle f,h \rangle$ and $k \neq 0$ so that
\[ w_{i-1} = \varphi(w) \, , \quad v_i = \varphi(v) \, , \quad v_{i+1} = \varphi f^k(v). \]
It follows by hypothesis (\ref{vwdistance}) that $\d(w_{i-1},v_i) = \d(w,v) \geq 2$, and once again $i(w_{i-1},v_i) \neq 0$.

Since $Z_i = \varphi(Z)$, hypothesis (\ref{vwprojections}) again shows
\begin{eqnarray*}
\d_{Z_i}(w_{i-1},v_{i+1}) & = & \d_{\varphi(Z)}(\varphi(w),\varphi f^k(v))\\
 & = & \d_Z(w,f^k(v))\\
 & >& 3M.\\
 \end{eqnarray*}

It follows that the sequence of special vertices for $\Phi_*(\gamma)$ satisfies the hypothesis of Proposition
\ref{geoconstruct}. Since $\Phi_*(\gamma)$ is obtained by concatenating geodesic segments between consecutive special
vertices, Proposition \ref{geoconstruct} completes the proof.
\end{proof}

We now turn our attention to arbitrary geodesics $\gamma$ of $\Ca$ (possibly infinite rays or biinfinite lines), and we would like to define a straightening for $\Phi(\gamma)$.
We do this first for a particular type of geodesic.
We say that $\gamma$ is {\em $h$--finite} if it contains no infinite $h$--geodesic ray.

Suppose now that $\gamma$ is $h$--finite.  Let $\gamma_1 \subset \gamma_2 \subset ... \subset \gamma$ be an exhaustion
by geodesic segments with the property that for each $i$, the first edge of $\gamma_i$ is an $h$--edge only if that
$h$--edge is the first edge of $\gamma$, and likewise, the last edge of $\gamma_i$ is an $h$--edge only if that
$h$--edge is the last edge of $\gamma$.  It follows that
\[
\Phi_*(\gamma_1) \subset \Phi_*(\gamma_2) \subset \cdots
\]
and we define $\Phi_*(\gamma)$ to be the union of these geodesic segments. Note that by construction,
$\Phi_*(\gamma)$ is (bi-)infinite if and only if $\gamma$ is. Thus we have

\begin{corollary} \label{hfinite}
If $\gamma$ is any $h$--finite geodesic in $\Ca$, then $\Phi_*(\gamma)$ is a geodesic in $\C(S)$, (bi-)infinite if and only if $\gamma$ is.
\end{corollary}

We can now also prove the first part of Theorem \ref{examplethrm}.

\begin{proposition} \label{injectivehom}
The canonical homomorphism $\langle f,h \rangle \to \Mod(S)$ is injective, and every element not conjugate to a power
of $h$ is pseudo-Anosov.
\end{proposition}

\begin{proof}
We suppose $\varphi$ is not conjugate to a power of $h$ and prove it is pseudo-Anosov.  This will prove the proposition.
Note that there is a biinfinite $h$--finite geodesic in $\Ca$ stabilized by $\varphi$, on which $\varphi$ acts by
translation. The straightening of its $\Phi$--image is a biinfinite geodesic in $\C(S)$ stabilized by $\varphi$. By
equivariance of $\Phi$, $\varphi$ acts by translation on this geodesic, and so it is pseudo-Anosov.
\end{proof}

Every point of $\partial \Ca$ is the endpoint of a unique geodesic ray beginning at $\1$. Denote the subset of
$\partial \Ca$ that are endpoints of $h$--finite geodesic rays beginning at $\1$ by $\partial^h \Ca$. From Corollary
\ref{hfinite}, we obtain a map from
\[
\partial \Phi_*: \partial^h \Ca \to \EL = \partial \C(S)
\]
sending the ideal endpoint of an $h$--finite geodesic ray $\gamma$ to the ideal endpoint of $\Phi_*(\gamma)$.
Note that this map is injective since any two distinct points $x,y \in \partial^h \Ca$ are the ideal endpoints of a
biinfinite geodesic $\gamma$.  Since $\gamma$ is clearly also $h$--finite, $\Phi_*(\gamma)$ is biinfinite with
distinct ideal endpoints $\partial \Phi_*(x)$ and $\partial \Phi_*(y)$, and so $\partial \Phi_*$ is injective.

Note that if $\gamma$ is any $h$--finite geodesic ray, by construction it contains infinitely many disjoint geodesic segments of length at least $6$
which are all $\langle f,h \rangle$--translates of $[v,f(v)]$.
In particular, there is a uniform bound on all projection coefficients for the endpoints of these segments, so by Corollary \ref{itworks}
we obtain the following.

\begin{proposition}
For every $x \in \partial^h \Ca$, $\partial \Phi_*(x)$ is uniquely ergodic.
\end{proposition}

We may therefore uniquely lift the map $\partial \Phi_*$ to a map (with the same name)
\[ \partial \Phi_*: \partial^h \Ca \to \PML \]
which is also injective.

Now we suppose $x \in \partial \Ca \setminus \partial^h \Ca$, and let $\gamma$ be a geodesic beginning at $\1$ ending
at $x$. Write $\gamma$ as the concatenation $\gamma = \gamma^\seg \gamma^\ray$, where $\gamma^\seg$ is a maximal
$h$--finite subgeodesic segment ending in an $f$--edge (which could be empty if $\gamma$ is an $h$--geodesic ray), and
$\gamma^\ray$ is an $h$--geodesic ray. Let $\varphi$ denote the terminal vertex of $\gamma^\seg$, and we define
$\Phi_*(\gamma)$ as the concatenated geodesic
\[
\Phi_*(\gamma) = \Phi_*(\gamma^\seg)[\varphi(v),\varphi(w)]
\]
which is indeed geodesic in $\C(S)$ by the same reasoning as above.

Now we exhaust the $h$--geodesic ray $\gamma^\ray$ by geodesic segments of the form $\gamma_k^\ray = \varphi([\1,h^k])$, for
$k \in {\mathbb Z}_+$ or $k \in {\mathbb Z}_-$, depending on whether the edges of the rays are positively or negatively oriented.
This provides an exhaustion of $\gamma$ by $h$--finite geodesic segments $\{ \gamma^\seg \gamma_k^\ray \}_{k \in {\mathbb Z}_{\pm}}$.
Furthermore, the associated geodesics are written as a concatenation
\[
\Phi_*(\gamma^\seg\gamma_k^\ray) = \Phi_*(\gamma)[\varphi(w),\varphi h^k(v)].
\]

We define $\partial \Phi_*(x)$ to be the $\varphi$--image of the stable lamination of $h$ if the ray is positively
oriented, and the $\varphi$-image of the unstable lamination of $h$ if it is negatively oriented. Equivalently, this is
the stable or unstable lamination of $\varphi h \varphi^{-1}$.

We observe that $\varphi(w)$, which is the terminal vertex of $\Phi_*(\gamma)$, is the unique curve disjoint from $\partial \Phi_*(x)$.
Thus, if $\partial \Phi_*(x) = \partial \Phi_*(y)$, and $\delta$ is the ray ending at $y$,
then the terminal vertex of $\Phi_*(\delta)$ is equal to that of $\Phi_*(\gamma)$.  By Lemma \ref{straightgeos},
we must have $\delta^\seg = \gamma^\seg$, and because the stable and unstable laminations of $h$
are distinct, it follows that $x = y$.  Thus we have proved that $\partial \Phi_*$ is injective on $\partial \Ca
\setminus \partial^h \Ca$.  Because these are non-filling laminations, while every lamination in $\partial
\Phi_*(\partial^h \Ca)$ is filling, this also proves

\begin{proposition}
$\partial \Phi_* : \partial \Ca \to \PML$ is injective.
\end{proposition}

All that remains is to prove the following.

\begin{proposition} \label{continuityprop}
$\partial \Phi_*: \partial \Ca \to \PML$ is continuous.
\end{proposition}

\begin{proof} We prove that $\partial \Phi_*$ is continuous at every $x \in \partial \Ca$.
The proof divides into two cases.\\

\noindent
{\bf Case 1.} $x \in \partial^h \Ca$.\\

We let $\{x_n \}_{n =1}^\infty \subset \partial \Ca$ with $x_n \to x$ as $n \to \infty$. By considering each situation
separately, we can assume that $\{x_n\}$ is completely contained in either $\partial^h \Ca$ or in the complement.
Let $\gamma_n$ and $\gamma$ be the geodesics beginning at $\1$ limiting on $x_n$ and $x$ for all $n$.\\

\noindent
{\bf Subcase 1.} $\{ x_n \} \subset \partial^h \Ca$.\\

Since $x_n \to x$, $\gamma_n$ converges to $\gamma$ uniformly on compact sets. Because $\Ca$ is a tree, it follows that
for any initial segment of $\gamma$, there is an initial segment of $\gamma_n$, for $n$ sufficiently large, which
agrees with this initial segment of $\gamma$.
Hence $\Phi_*(\gamma_n)$ converges uniformly on compact sets to $\Phi_*(\gamma)$, and $\partial \Phi_*(x_n) \to \partial \Phi_*(x)$, as required.\\

\noindent
{\bf Subcase 2.} $\{ x_n \}_{n=1}^\infty \subset \partial \Ca \setminus \partial^h \Ca$.\\

Since each $\gamma_n$ is an $h$--infinite geodesic, $\Phi_*(\gamma_n)$ is a finite geodesic whose terminal vertex we
denote $w_n$ (which is disjoint from $\partial \Phi_*(x_n)$).

Because $x_n \to x$, we again see that  $\Phi_*(\gamma_n)$ converges on compact sets to $\Phi_*(\gamma)$.
Since $\Phi_*(\gamma)$ is an infinite geodesic ray, it follows that the endpoints $w_n$ of $\Phi_*(\gamma_n)$ converge to $\partial \Phi_*(x)$.
Because $\partial \Phi_*(x_n)$ is disjoint from $w_n$, every accumulation point of $\{ \partial \Phi_*(x_n) \}$ has intersection number zero with $\partial \Phi_*(x)$.
Finally, the fact that $\partial \Phi_*(x)$ is uniquely ergodic implies
\[ \lim_{n \to \infty} \partial \Phi_*(x_n) = \partial \Phi_*(x) \]
as required.\\

\noindent
{\bf Case 2.} $x \in \partial \Ca \setminus \partial^h \Ca$.\\

Again, suppose $\{ x_n \} \subset \partial \Ca$ is a sequence converging to $x$, and let $\gamma_n$ and $\gamma$ be
geodesic rays limiting on $x_n$ and $x$, respectively, for all $n$. Since $\gamma$ is $h$--infinite, $\Phi_*(\gamma)$
is finite and we let $\varphi$ denote the terminal vertex of $\gamma^\seg$ (notation as above) so that $\varphi(w)$ is
the terminal vertex of $\Phi_*(\gamma)$ (which is the unique curve disjoint $\partial \Phi_*(x)$).

As above, since $\Ca$ is a tree and $x_n \to x$, it follows that any initial segment of $\gamma$ is equal to some initial segment of $\gamma_n$
for all sufficiently large $n$.
By throwing away finitely many initial terms in the sequence $\gamma_n$ (which we can do without loss of generality) we decompose each
$\gamma_n$ as a concatenation
\[ \gamma_n = \gamma^\seg \gamma_n^h \gamma_n^\infty \]
where $\gamma_n^h$ is the segment of $\gamma^\ray$ such that $\gamma^\seg \gamma_n^h$ is the largest segment of $\gamma_n$ that agrees with an
initial segment of $\gamma$.
The ray $\gamma_n^\infty$ is then $\overline{\gamma_n - (\gamma^\seg \gamma_n^h)}$, and we note that its initial edge is an $f$--edge by maximality
of $\gamma^\seg \gamma_n^h$.

We can then express the geodesics in $\C(S)$ associated to these $\gamma_n$ as
\[
\Phi_*(\gamma_n) = \Phi_*(\gamma)[\varphi(w),\varphi h^{k(n)}(v)] \Phi_*(\gamma_n^\infty).
\]
Since the $\gamma_n$ agree with $\gamma$ on longer and longer segments, it follows that $k(n) \to \infty$ or $k(n) \to
-\infty$ as $n \to \infty$, depending on the orientation of $\gamma^\ray$. We assume $k(n) \to \infty$, the
other case being similar.

As $n \to \infty$, $\varphi h^{k(n)}(v) = \varphi h^{k(n)} \varphi^{-1}(\varphi(v))$ tends to $\partial \Phi_*(x)$ (the stable lamination of
$\varphi h \varphi^{-1}$). Setting $Y_0 = S \setminus \varphi(w) = \varphi(Y)$, which is the supporting subsurface of $\partial \Phi_*(x)$, Theorem
\ref{bgithrm} implies
\[
\d_{Y_0}(\varphi h^{k(n)} (v),\partial \Phi_*(x_n)) \leq \diam_{Y_0}(\Phi_*(\gamma_n^\infty)) \leq M.
\]
Since $\pi_{Y_0}(\varphi h^{k(n)}(v))$ is tending to $|\partial \Phi_*(x)|$ in $\A(Y_0) \cup \EL(Y_0)$, it follows that any accumulation point of $\partial \Phi_*(x_n)$ in $\PML(S)$ must have zero intersection number with $\partial \Phi_*(x)$.
Since $\partial \Phi_*(x)$ is uniquely ergodic (though not filling) we see that any limit of $\partial \Phi_*(x_n)$ is a point in the projective $1$--simplex
of measures supported on $\varphi(w) \cup |\partial \Phi_*(x)|$.

We suppose $\mu$ is any limit of $\partial \Phi_*(x_n)$ and show that the support is $|\partial \Phi_*(x)|$, which will
complete the proof. Replace $\partial \Phi_*(x_n)$ by a subsequence that converges to $\mu$, and we further assume (by
possibly replacing with a smaller subsequence) that the Hausdorff limit of $|\partial \Phi_*(x_n)|$ exists. Note that
if $\mu$ had some non-trivial transverse measure on $\varphi(w)$, then the Hausdorff limit of $|\partial \Phi_*(x_n)|$
would contain $\varphi(w)$. If this were true, then it must be the case that
\[ \d_{\varphi(w)}(\varphi(v),\partial \Phi_*(x_n)) \to \infty \]
as $n \to \infty$. However, hypothesis (\ref{htranslation}) implies
\[ \d_{\varphi(w)}(\varphi(v),\varphi h^{k(n)}(v)) = \d_{w}(v,h^{k(n)}(v)) \leq 2 \]
and so combined with the triangle inequality and Theorem \ref{bgithrm} we obtain
\begin{eqnarray*}
\d_{\varphi(w)}(\varphi(v),\partial \Phi_*(x_n)) & \leq & \d_{\varphi(w)}(\varphi(v),\varphi h^{k(n)}(v)) + \d_{\varphi(w)}(\varphi h^{k(n)}(v),\partial \Phi_*(x_n))\\
 & \leq & 2 + \diam_{\varphi(w)}(\Phi_*(\gamma_n^\infty))\\
 & \leq & 2 + M.
\end{eqnarray*}
Therefore, $\mu$ has no measure on $\varphi(w)$, and hence is supported on $|\partial \Phi_*(x)|$, completing the proof.
\end{proof}

\subsection{Constructing $f$ and $h$} \label{theyexist}

Let us now explain how to find $f$ and $h$ satisfying all the hypotheses.

Finding $h$ satisfying hypothesis (\ref{hreduce}) is easy, and the first part of (\ref{htranslation}) is gotten by replacing $h$ by
any sufficiently large power.
To guarantee that $h$ satisfies the second part of (\ref{htranslation}) first replace $h$ by a power that leaves invariant each boundary leaf of the
stable lamination $|\lambda|$ for $h$.
The component $U$ of the path metric completion of $S \setminus |\lambda|$ containing $w$ is a union of two {\em crowns} along $w$; see
\cite{cassonbleiler}.
Let $\ell \subset U$ denote a biinfinite geodesic passing through $w$ exactly once and running from the cusp of one crown to the cusp of another.
Denoting the Dehn twist in $w$ by $T_w$, we replace $h$ by $T_w^k h$, for an appropriate $k$ so that $h(\ell) = \ell$, so that hypothesis (\ref{htranslation}) is satisfied.

Hypothesis (\ref{fpseudo}) is easily arranged by assuming the stable and unstable laminations for a pseudo-Anosov mapping
class $f$ have a big projection coefficient to the complement of a nonseparating curve $v$.  For then, after possibly
replacing $f$ by a sufficiently large power and connecting $v$ to $f(v)$ by a geodesic, taking the $f$--orbit gives a biinfinite geodesic
$\g$ as required (this is a geodesic by Proposition \ref{geoconstruct}). Replacing $f$ by a larger power, we also
guarantee that hypothesis (\ref{ftranslation}) is satisfied.

Replacing $f$ with a conjugate by a sufficiently large power of an independent pseudo-Anosov mapping class $\varphi \in \Mod(S)$
will ensure $\d(v', w) \geq 2$ for every $v' \in \g$.
We further replace $f$ by its conjugate by a large power of any element $\varphi \in \Mod(S)$ which is pseudo-Anosov on the complement of $v$.
Taking this power sufficiently large, we guarantee that for any $v' \in \g$ we have $\d_Z(v',w) > 3M$, and hence by Theorem \ref{bgithrm} a geodesic from $w$ to $v'$
passes through $v$.
In particular, we have $\d(w,v') \geq \d(w,v) \geq 2$, guaranteeing hypothesis (\ref{vwdistance}), and since $v' = f^k(v)$ is a vertex of $\g$,
we have also arranged hypothesis (\ref{vwprojections}).


\bibliographystyle{plain}
\bibliography{triplebib}

\bigskip

\noindent Department of Mathematics, Brown University, Providence, RI 02912 \newline \noindent
\texttt{rkent@math.brown.edu}

\bigskip

\noindent Department of Mathematics, University of Illinois, Urbana-Champaign, IL 61801 \newline \noindent  \texttt{clein@math.uiuc.edu}

\end{document}